\documentclass[12pt]{article}
\usepackage{etex,latexsym,amsmath,amssymb,amsthm}
\usepackage{graphtex}
\usepackage{multicol,epsfig,array}
\usepackage{cancel,multirow}

\title{\bf Revisiting a Number-Theoretic Puzzle: \\ The Census-Taker Problem}
\author{{\bf I.J.L. Garces} \\ {\bf Mark L. Loyola} \\ Department of Mathematics \\ Ateneo de Manila University}
\date{27 July 2010}

\begin{document}
\maketitle

\newtheorem{example}{Example}
\newtheorem{exercise}[example]{Exercise}
\newtheorem{theorem}[example]{Theorem}
\newtheorem{lemma}[example]{Lemma}
\newtheorem{corollary}[example]{Corollary}

\def\ds{\displaystyle}

\section{Introduction}
A census taker knocks on a door. A mother answers.

\medskip
The census taker says, ``I need to know the number of children you have, and their ages.'' The woman responds in puzzle-ese, ``I have three daughters, the product of their ages is 36, and the sum of their ages is equal to the house number next door.''

\medskip
The census taker, who never wastes questions, computes for a while and then asks, ``Does your oldest daughter love dogs?'' The mother answers affirmatively. The census taker says, ``Thank you. I now know the ages.''

\medskip
What are the ages of the children?

\medskip
With seemingly insufficient information, such number-theoretic puzzles belong to what is now known as the {\it census-taker problem}.

\medskip
To someone encountering the problem for the first time, the information transpired during the conversation is not sufficient to determine the ages. Aside from this initial reaction, many questions arise. How can one determine the ages from the given product and the house number next door? What does the oldest daughter's loving of dogs have to do with the problem? To solve the problem, one simply puts his feet on the shoes of the census taker. The readers do not know the house number next door, but the census taker does. Moreover, once one lists down the eight triples of positive integers whose product is $36$, the answers to these questions will become clearer. This list is presented in Table \ref{listoftriplescensustaker36}, where each triple consists of the ages of the daughters. Looking at the table, there are two triples whose sums are equal to $13$. Because the census taker never wastes a question, this common sum is the reason why he asks another question, and this sum
 is also the house number next door. The affirmative answer of the mother confirms that she has an oldest daughter, which nails down the triple $\{9,2,2\}$ to give the ages of the girls.

\begin{table}[h]
$$\begin{tabular}{|c|c||c|c|} \hline
{\bf Triple} & {\bf Sum} & {\bf Triple} & {\bf Sum} \\ \hline\hline
$\{36,1,1\}$ & $38$ & $\{9,2,2\}$ & $13$ \\
$\{18,2,1\}$ & $21$ & $\{6,6,1\}$ & $13$ \\
$\{12,3,1\}$ & $16$ & $\{6,3,2\}$ & $11$ \\
$\{9,4,1\}$ & $14$ & $\{4,3,3\}$ & $10$ \\ \hline
\end{tabular}$$
\caption{Triples whose product is $36$.}\label{listoftriplescensustaker36}
\end{table}

The last question of the census taker is important. Its answer decides which of the triples $\{9,2,2\}$ and $\{6,6,1\}$ gives the ages of the daughters. What if the mother's answer to the last question is negative? Does it point to the triple $\{6,6,1\}$? Not necessarily. It depends on how the mother phrases her negative answer. If she simply says ``No,'' this means that she has an oldest daughter who does not love dogs, which points again to the triple $\{9,2,2\}$. However, if the mother says something like ``No, but my youngest daughter does,'' then this answer confirms that she does have a youngest daughter, which points to the triple $\{6,6,1\}$. Asking the last question then only tells us that the census taker is choosing between two triples, but the mother's appropriate answer settles the dilemma.

\medskip
According to \cite{greenbelt}, the census-taker problem originated during World War II. In literature, it has appeared in different stories in a puzzle-ese style. As a result, many number-theoretic problems with seemingly insufficient information have been proposed in some journals. See \cite{problem1} and \cite{problem2}.

\medskip
The current work revisits the results in \cite{meyers&see}, and presents the census-taker problem as a motivation to introduce the beautiful theory of numbers.

\section{The Census-Taker Numbers}
Table \ref{listoftriplescensustaker36} depicts a special property of the number $36$. Among the distinct unordered triples of positive integers whose product is $36$, there are exactly two triples with equal sums. With such special property, the number $36$ is called a {\it census-taker number}, and $36$ is the smallest such positive integer.

\medskip
More formally, a positive integer $N$ is said to be a {\it census-taker number} (CTN) if there is exactly one pair $\{\{a,b,c\},\{d,e,f\}\}$ of different (unordered) triples of positive integers such that
\begin{equation}\label{productcondition}
abc=def=N
\end{equation}
and
\begin{equation}\label{sumcondition}
a+b+c=d+e+f.
\end{equation}
We shall refer equation (\ref{productcondition}) as the {\it product condition}, equation (\ref{sumcondition}) as the {\it sum condition}, the common sum as the {\it magic sum}, and the triples $\{a,b,c\}$ and $\{d,e,f\}$ as the {\it mysterious triples} of the CTN.

\medskip
It can be easily checked that there are no one-digit CTNs. Moreover, Table \ref{twodigitcensustakernumbers} lists down all the two-digit CTNs. For the curious minds, there are twenty-nine three-digit CTNs.

\begin{table}[h]
$$\begin{tabular}{|c|c|c|c|} \hline
{\bf CTN} & {\bf Magic Sum} & {$\boldsymbol{\{a,b,c\}}$} & {$\boldsymbol{\{d,e,f\}}$} \\ \hline\hline
$36$ & $13$ & $\{9,2,2\}$ & $\{6,6,1\}$ \\
$40$ & $14$ & $\{10,2,2\}$ & $\{8,5,1\}$ \\
$72$ & $14$ & $\{8,3,3\}$ & $\{6,6,2\}$ \\
$96$ & $21$ & $\{16,3,2\}$ & $\{12,8,1\}$ \\ \hline
\end{tabular}$$
\caption{The two-digit CTNs, and their magic sums and mysterious triples.}\label{twodigitcensustakernumbers}
\end{table}

\begin{figure}[h]
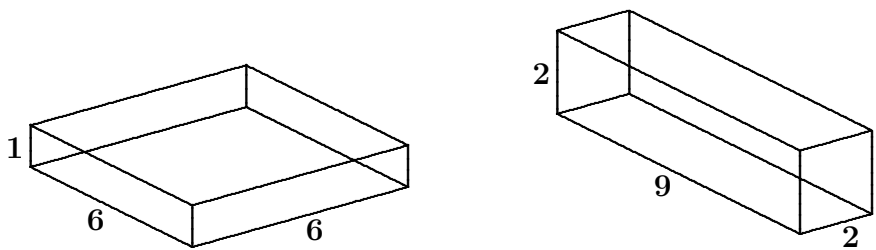

$$\pic
\xunit17pt \yunit17pt \zunit17pt

\PlotSize4


\LINE (0,0,0) (6,0,0) (6,6,0) (0,6,0) (0,0,0)
\LINE (0,0,1) (6,0,1) (6,6,1) (0,6,1) (0,0,1)
\LINE (0,0,0) (0,0,1)
\LINE (6,0,0) (6,0,1)
\LINE (6,6,0) (6,6,1)
\LINE (0,6,0) (0,6,1)

\ALIGN[r] ($\mathbf{1}$) (6.2,0,0.5)
\ALIGN[ur] ($\mathbf{6}$) (6.2,3,0)
\ALIGN[ul] ($\mathbf{6}$) (3,6.2,0)
\cip \qquad\qquad
\pic
\xunit17pt \yunit17pt \zunit17pt

\PlotSize4


\LINE (0,0,0) (2,0,0) (2,9,0) (0,9,0) (0,0,0)
\LINE (0,0,2) (2,0,2) (2,9,2) (0,9,2) (0,0,2)
\LINE (0,0,0) (0,0,2)
\LINE (2,0,0) (2,0,2)
\LINE (2,9,0) (2,9,2)
\LINE (0,9,0) (0,9,2)

\ALIGN[r] ($\mathbf{2}$) (2.2,0,1)
\ALIGN[ur] ($\mathbf{9}$) (2.2,4.5,0)
\ALIGN[ul] ($\mathbf{2}$) (1,9.2,0)
\cip$$
\caption{The ``mysterious'' rectangular boxes with volume equal to $36$.}\label{boxeswithcensustakervolume}
\end{figure}

There is an interesting geometric interpretation of the census-taker number. Every CTN produces a unique pair of ``mysterious'' rectangular boxes with special properties. Each of the mysterious triples constructs a corresponding mysterious rectangular box whose dimensions are the elements of the corresponding triple. The product condition implies that these boxes have the same volume (which is the given CTN), while the sum condition implies that the boxes have the same perimeter (which is four times the magic sum). The volume being a census-taker number denies the existence of a third box of integral dimensions with the same volume and perimeter. Figure \ref{boxeswithcensustakervolume} shows the two mysterious rectangular boxes whose dimensions correspond to the mysterious triples $\{9,2,2\}$ and $\{6,6,1\}$ of the census-taker number $36$.

To show that a number is a census-taker, it is not enough to find two triples that satisfy the product and the sum conditions. For example, consider $N=144$. Note that
$$144 = 9 \cdot 4 \cdot 4 = 8 \cdot 6 \cdot 3 \qquad\mbox{and}\qquad 9+4+4=8+6+3,$$
and so the triples $\{9,4,4\}$ and $\{8,6,3\}$ satisfy the product and the sum conditions. However, $N=144$ is not a census-taker number because another pair of triples satisfies the product and the sum conditions, and these are $\{12,4,3\}$ and $\{9,8,2\}$. For a positive integer to be a census-taker number, there must be exactly two triples that satisfy the product and the sum conditions.

\section{Properties of CTNs}
Table \ref{twodigitcensustakernumbers} listed down the four two-digit CTNs, and it was mentioned that there are twenty-nine three-digit CTNs. For smaller integers like $36$, the determination of the CTNs can be done by enumerating and then checking every possible triple. However, for larger integers, the impracticality of this procedure demands a more elegant way to determine the CTNs. Some nice properties of these numbers help us eliminate most of those integers that are not CTNs in a faster way. As concrete illustration, we apply this sieving technique to positive integers up to $100$.

\medskip
Let $\{a,b,c\}$ and $\{d,e,f\}$ be the mysterious triples associated with the census-taker number $N$. Since the order of the triples and the order of the elements in each triple are not important, we assume from now on that
$$a \ge b \ge c, \quad d \ge e \ge f, \quad\mbox{and}\quad c \ge f.$$
This assumption will be further restricted after the following theorem.

\begin{theorem}\label{triplesnocommonentry}
The mysterious triples of a CTN cannot have a common entry.
\end{theorem}

{\it Proof.} Let $\{a,b,c\}$ and $\{d,e,f\}$ be the mysterious triples of the CTN. Without loss of generality, we suppose that $a=d$. By the sum and the product conditions, we have
$$f-c=b-e \qquad\mbox{and}\qquad bc=ef.$$
Using the identities $bc=ec+c(b-e)$ and $ef=ec+e(f-c)$, we get
$$c(b-e)=e(b-e),$$
which implies that $b=e$ or $c=e$. The former implies that $c=f$, while the latter implies that $b=c=e=f$. In both cases, we conclude that $\{a,b,c\}=\{d,e,f\}$, a contradiction. \hfill {\sc q.e.d.}

\medskip
With Theorem \ref{triplesnocommonentry}, we further assume that $c > f$. More than this assumption, Theorem \ref{triplesnocommonentry} eliminates a good number positive integers as candidates for CTNs. A prime cannot be a CTN because all triples of positive integers whose product is the given prime have the integer $1$ as a common element. In fact, this quick consequence of Theorem \ref{triplesnocommonentry} also eliminates any integer of the form $pq$, where $p$ and $q$ are (not necessarily distinct) primes, from the list of CTNs. This is so because the only possible triples for these integers $\{pq,1,1\}$ and $\{p,q,1\}$, and these triples have a common element $1$.

\begin{theorem}\label{CTNnotprimepower}
A CTN is not a power of a prime.
\end{theorem}

{\it Proof.} Suppose, on the contrary, that $N=p^n$ is a CTN, where $p$ is prime and $n$ is a positive integer. Since any divisor of $N$ is also a power of $p$, we may assume that $\{p^r,p^s,p^t\}$ and $\{p^u,p^v,p^w\}$ are the mysterious triples of $N$, where $n=r+s+t=u+v+w$, $r \ge s \ge t$, $u \ge v \ge w$, and $t > w$. From the sum condition, we have
$$p^{t-w}(p^{r-t}+p^{s-t}+1)=p^{u-w}+p^{v-w}+1. \eqno{(\star)}$$
Note that the left-hand side of equation $(\star)$ is divisible by $p$. However, for $p \ne 3$, the right-hand side of $(\star)$ is not divisible by $p$. Finally, for $p=3$, the right-hand side is divisible by $p$ if and only if $u=v=w$, which implies that the right-hand is exactly equal to $3$, while the left-hand side is at least $9$. Both scenarios lead to a contradiction. \hfill {\sc q.e.d.}

\begin{theorem}\label{CTNatleastfourprimes}
Every CTN is a product of at least four (not necessary distinct) primes.
\end{theorem}

{\it Proof.} By Theorem \ref{triplesnocommonentry}, it suffices to prove that a CTN is not a product of exactly three primes.

Suppose that $N=pqr$ is a CTN, where $p$, $q$, and $r$ are primes with $p \ge q \ge r$. Then there are only five possible triples for $N$, and these are $\{pqr,1,1\}$, $\{pq,r,1\}$, $\{pr,q,1\}$, either $\{qr,p,1\}$ or $\{p,qr,1\}$, and $\{p,q,r\}$. By Theorem \ref{triplesnocommonentry}, the only possible pair of mysterious triples is $\{pqr,1,1\}$ and $\{p,q,r\}$. However, since $p+q+r < pqr+2$, this pair does not satisfy the sum condition. \hfill {\sc q.e.d.}

\medskip
Table \ref{possibleCTNs} shows the remaining integers from $1$ to $100$ after applying Theorems \ref{triplesnocommonentry}, \ref{CTNnotprimepower}, and \ref{CTNatleastfourprimes}, where the integers are crossed out with $/$, $//$, and ---, respectively.

\begin{table}[h]
$$\begin{array}{cccccccccc}
\cancel{1} & \cancel{11} & \cancel{21} & \cancel{31} & \cancel{41} & \cancel{51} & \cancel{61} & \cancel{71} & \rlap{//}{81} & \cancel{91} \\
\cancel{2} & \rlap{---}{12} & \cancel{22} & \rlap{//}{32} & \rlap{---}{42} & \rlap{---}{52} & \cancel{62} & 72 & \cancel{82} & \rlap{---}{92} \\
\cancel{3} & \cancel{13} & \cancel{23} & \cancel{33} & \cancel{43} & \cancel{53} & \rlap{---}{63} & \cancel{73} & \cancel{83} & \cancel{93} \\
\cancel{4} & \cancel{14} & 24 & \cancel{34} & \rlap{---}{44} & 54 & \rlap{//}{64} & \cancel{74} & 84 & \cancel{94} \\
\cancel{5} & \cancel{15} & \cancel{25} & \cancel{35} & \rlap{---}{45} & \cancel{55} & \cancel{65} & \rlap{---}{75} & \cancel{85} & \cancel{95} \\
\cancel{6} & \rlap{//}{16} & \cancel{26} & 36 & \cancel{46} & 56 & \rlap{---}{66} & \rlap{---}\rlap{---}{76} & \cancel{86} & 96 \\
\cancel{7} & \cancel{17} & \rlap{//}{27} & \cancel{37} & \cancel{47} & \cancel{57} & \cancel{67} & \cancel{77} & \cancel{87} & \cancel{97} \\
\rlap{\!/\!/}{8} & \rlap{---}{18} & \rlap{---}{28} & \cancel{38} & 48 & \cancel{58} & \rlap{---}{68} & \rlap{---}{78} & 88 & \rlap{---}{98} \\
\cancel{9} & \cancel{19} & \cancel{29} & \cancel{39} & \cancel{49} & \cancel{59} & \cancel{69} & \cancel{79} & \cancel{89} & \rlap{---}{99} \\
\cancel{10} & \rlap{---}{20} & \rlap{---}{30} & 40 & \rlap{---}{50} & 60 & \rlap{---}{70} & 80 & 90 & 100
\end{array}$$
\caption{Applying Theorems \ref{triplesnocommonentry}, \ref{CTNnotprimepower}, and \ref{CTNatleastfourprimes} to integers from $1$ to $100$.}\label{possibleCTNs}
\end{table}

In Table \ref{possibleCTNs}, fourteen two-digit integers remain after the sieving techniques due to the theorems. As claimed in Table \ref{twodigitcensustakernumbers}, ten of these integers are not CTNs. Albeit there are advanced properties of the CTNs that are capable of eliminating some of these from the list, we do not mention them anymore in this paper, and we remain with the elementary treatment of these interesting numbers. To check whether or not the remaining fourteen integers in Table \ref{possibleCTNs} are CTNs, we simply list down all the triples for each remaining integer, and verify the definition. As examples, Table \ref{listoftriples56} shows that $N=56$ is not a CTN by simply looking at the different sums produced by the triples, and Table \ref{listoftriples96} proves that $N=96$ is a CTN.

\begin{table}[h]
$$\begin{tabular}{|c|c||c|c|} \hline
{\bf Triple} & {\bf Sum} & {\bf Triple} & {\bf Sum} \\ \hline\hline
$\{56,1,1\}$ & $58$ & $\{14,2,2\}$ & $18$ \\
$\{28,2,1\}$ & $31$ & $\{8,7,1\}$ & $16$ \\
$\{14,4,1\}$ & $19$ & $\{7,4,2\}$ & $13$ \\ \hline
\end{tabular}$$
\caption{Triples whose product is $56$.}\label{listoftriples56}
\end{table}
\begin{table}[h]
$$\begin{tabular}{|c|c||c|c|} \hline
{\bf Triple} & {\bf Sum} & {\bf Triple} & {\bf Sum} \\ \hline\hline
$\{96,1,1\}$ & $98$ & $\{16,3,2\}$ & $21$ \\
$\{48,2,1\}$ & $51$ & $\{12,8,1\}$ & $21$ \\
$\{32,3,1\}$ & $36$ & $\{12,4,2\}$ & $18$ \\

$\{24,4,1\}$ & $29$ & $\{8,6,2\}$ & $16$ \\
$\{24,2,2\}$ & $28$ & $\{8,4,3\}$ & $15$ \\
$\{16,6,1\}$ & $23$ & $\{6,4,4\}$ & $14$\\ \hline
\end{tabular}$$
\caption{Triples whose product is $96$.}\label{listoftriples96}
\end{table}

\section{How many CTNs are there?}
It is not known whether or not there are infinitely many census-taker numbers. However, if we assume that there are infinitely many primes $p$ such that $2p-1$ is also prime, we claim that there are also infinitely many CTNs.

\medskip
We need the following inequalities to prove our claim about the infinitude of CTNs.

\begin{lemma}\label{lemma1}
For all real numbers $x \ne 1$, the following inequalities hold:
\begin{enumerate}
\item[{(i)}] $4x^4-4x^3+x^2+2 > 4x^3-4x^2+2x+1$
\item[{(ii)}] $5x^2-4x+2 > 4x^2-2x+1$
\item[{(iii)}] $4x^2-2x+1 > 2x^2+2x-1$
\item[{(iv)}] $2x^2+2x-1 > x^2+4x-2$
\end{enumerate}
\end{lemma}

{\it Proof.} To present the proof technique to establish these inequalities, we completely prove (i).

Observe that
$$4x^4-4x^3+x^2+2 = 4x^3-4x^2+2x+1 + (4x^2+1)(x-1)^2.$$
Since $(4x^2+1)(x-1)^2>0$ for all $x \ne 1$, inequality (i) follows.

For the other inequalities, we simply observe the following equations:
\begin{align*}
5x^2-4x+2 &= 4x^2-2x+1 + (x-1)^2 \\
4x^2-2x+1 &= 2x^2+2x-1 + 2(x-1)^2 \\
2x^2+2x-1 &= x^2+4x-2 + (x-1)^2.
\end{align*}
Since $(x-1)^2>0$ for all $x \ne 1$, the other inequalities follow. \hfill {\sc q.e.d.}

\medskip
With slight modifications, the technique used in the proof of the previous lemma can also be utilized to prove the following lemmas.

\begin{lemma}\label{lemma2}
For all real numbers $x>1$, the following inequality holds:
$$2x^3-x^2+2x > 5x^2-4x+2.$$
\end{lemma}

\begin{lemma}\label{lemma3}
For all real numbers $x$ that satisfy $x > -\frac{1}{2}$ and $x \ne 1$, the following inequality holds:
$$4x^3-4x^2+2x+1 > 2x^3-x^2+2x.$$
\end{lemma}

We are now ready to prove a theorem that will lead to our claim.

\begin{theorem}\label{aprimeimplyingCTN}
Let $p$ be prime such that $2p-1$ is also prime. Then $N=p^2(2p-1)^2$ is a CTN.
\end{theorem}

{\it Proof.} Table \ref{listoftriplesp^2(2p-1)^2} presents the eight possible triples of $N$.
\begin{table}[t]
$$\begin{tabular}{|c|c||c|c|} \hline
{\bf Triple} & {\bf Sum} & {\bf Triple} & {\bf Sum} \\ \hline\hline
$\{p^2(2p-1)^2,1,1\}$ & $4p^4-4p^3+p^2+2$ & $\{(2p-1)^2,p,p\}$ & $4p^2-2p+1$ \\
$\{p(2p-1)^2,p,1\}$ & $4p^3-4p^2+2p+1$ & $\{p(2p-1),p(2p-1),1\}$ & $4p^2-2p+1$ \\
$\{p^2(2p-1),2p-1,1\}$ & $2p^3-p^2+2p$ & $\{p(2p-1),2p-1,p\}$ & $2p^2+2p-1$ \\
$\{(2p-1)^2,p^2,1\}$ & $5p^2-4p+2$ & $\{p^2,2p-1,2p-1\}$ & $p^2+4p-2$ \\ \hline
\end{tabular}$$
\caption{Triples whose product is $p^2(2p-1)^2$, where $p$ and $2p-1$ are primes.}\label{listoftriplesp^2(2p-1)^2}
\end{table}
Using Lemmas \ref{lemma1}, \ref{lemma2}, and \ref{lemma3}, the sums can be arranged from highest to lowest in the following inequalities:
\begin{align*}
4p^4-4p^3+p^2+2 &> 4p^3-4p^2+2p+1 \\
&> 2p^3-p^2+2p \\
&> 5p^2-4p+2 \\
&> 4p^2-2p+1 \\
&> 2p^2+2p-1 \\
&> p^2+4p-2.
\end{align*}
It follows that $N=p^2(2p-1)^2$ is a CTN. \hfill {\sc q.e.d.}

\begin{corollary}\label{infinitudeofCTN}
If there are infinitely many primes $p$ such that $2p-1$ is also prime, then there are also infinitely many CTNs.
\end{corollary}

As a serious remark to the infinitude of CTNs established in the above corollary, it is not yet known, however, whether or not there are indeed infinitely many such primes $p$. Nevertheless, because there are CTNs that are not of the form $p^2(2p-1)^2$, Corollary \ref{infinitudeofCTN} is not the only direction to prove or disprove the infinitude of CTNs.

\section{Concluding Remarks}
A related problem was studied in \cite{kelly}, where it was shown that, for every integer $M > 18$, there exist triples with sum $M$ and equal (unspecified) products. Though this problem can be thought of as a dual of the census-taker problem, it does not consider the uniqueness condition that there should be exactly two triples with equal products. In fact, this uniqueness condition is what makes the census-taker problem more interesting and difficult.

\medskip
For the curious minds, we mention again that there are no one-digit, $4$ two-digit, and $29$ three-digit CTNs. Adding to these lists, with $9996$ as the largest, there are $277$ four-digit CTNs.

\medskip
With further twists and turns, the census-taker problem still has a long journey to go.

\bigskip\noindent
{\bf Acknowledgments.} The authors thank Giovanni Mazzarello for generously sharing his programming skills by providing the list of all CTNs, with magic sums and mysterious triples, up to 10000. Moreover, they also thank the anonymous referees for the constructive reviews, especially the referee who suggested to include a separate section about the number of CTNs.

\end{document}